\documentclass[11pt,reqno]{amsart}

\newcommand{\T}{{\mathbf T}^m}

\newcommand{\K}{{\mathbf K}^\crit}

\newcommand{\RR}{{\mathbb R}}

\newcommand{\ZZ}{{\mathbb Z}}

\newcommand{\NN}{{\mathbb N}}

\newcommand{\R}{{\mathbb R}}
\newcommand{\C}{{\mathbb C}}

\newcommand{\Z}{{\mathbb Z}}

\newcommand{\N}{{\mathbb N}}

\newcommand{\diag}{{\operatorname{diag}}}

\newcommand{\half}{{\frac{1}{2}}}

\newcommand{\be}{\beta}
\newcommand{\ga}{\gamma}

\def    \half   {{\frac{1}{2}}}

\def    \Z  {{\mathbb Z}}
\def    \R  {{\mathbb R}}
\def    \C  {{\mathbb C}}

 \def   \half   {{\frac{1}{2}}}

\newtheorem{theo}{{\sc Theorem}}[section]
\newtheorem{cor}[theo]{{\sc Corollary}}
\newtheorem{conj}[theo]{{\sc Conjecture}}
\newtheorem{remark}[theo]{{\sc Remark}}
\newtheorem{lem}[theo]{{\sc Lemma}}

\def\h#1{\hbox{#1}}

\def\o{\omega}
\def\O{\Omega}

\def\K{K\"ahler }

\def\ra{\rightarrow}
\def\a{\alpha}
\def\th{\theta}
\def\vp{\varphi}
\def\pa{\partial}
\def\isom{\cong}

\def\w{\wedge}
\def\i{\sqrt{-1}}
\def\text{\textstyle}

\def\ra{\rightarrow}

\def\isom{\cong}

\chardef\inodot="10
\def\q{\quad}
\def\qq{\qquad}
\def\all{\forall}
\def\opcit{\underbar{\phantom{aaaaa}}}

\def\dis{\displaystyle}

\def\del{\partial}

\def\frac#1#2{{{#1}\over{#2}}}
\def\ov{\over}
\def\dis{\displaystyle}
\def\h{\hbox}

\def\ppx{\frac{\partial}{\partial x}}
\def\ppt{\frac{\partial}{\partial t}}

\def\N{\nabla}
\def\K{K\"ahler }
\def\Sch{Schr\"odinger }
\def\bfa{{\bf a}}

\def\Nt{\nabla_t}
\def\Nx{\nabla_x}
\def\i{\sqrt{-1}}
\def\G{\Gamma}
\def\T#1#2#3{\G^{#1}_{#2#3}}

\def\be{\beta}
\def\a{\alpha}
\def\b#1{\bar{#1}}

\def\ustarTN{u^\star TN}

\def\TstarN{T^\star N}
\def\TstN{T^\star N}
\def\ustarh{u^\star h}
\def\dVMg{dV_{M,g}}
\def\dVNh{dV_{N,h}}
\def\pa{\partial}
\def\dudt{\frac{\partial u}{\pa t}}

\def\ddto#1{\frac{\partial #1}{\pa t}\big|_{0}}
\def\ddtBigo#1{\frac{\partial #1}{\pa t}\Big|_{0}}
\def\ddto#1{\frac{\partial #1}{\pa t}\big|_{0}}
\def\ustar{u^\star}
\def\tr{\h{tr}}

\def\bfP{{\bf P}}

\def\th{\theta}

\def\n#1#2{||#1||_{#2}}

\def\N{\nabla}

\def\bfPhi{{\bf \Phi}} \def\bfS{{\bf S}}   \def\bfT{{\bf T}}
\def\bfPhitilde{{\bf \tilde\Phi}}
\def\strutdepth{\dp\strutbox}
\def\specialstar{\vtop to \strutdepth{
    \baselineskip\strutdepth
    \vss\llap{$\star$\ \ \ \ \ \ \ \ \  }\null}}
\def\marginalstar{\strut\vadjust{\kern-\strutdepth\specialstar}}

\def\ga{\gamma}
\def\bfQ{{\bf Q}}
\def\bfP{{\bf P}}
\def\bfPhi{{\bf \Phi}} \def\bfS{{\bf S}}   \def\bfT{{\bf T}}

\def\parr#1#2{\frac{\partial#1}{\partial {#2}}}
\def\isom{\cong}

\def\pt{\h{pt}}

\title[One-dimensional Schr\"odinger map flow]
{On the global well-posedness of the one-dimensional Schr\"odinger map flow}

\author{Igor Rodnianski }
\author{Yanir A. Rubinstein }
\author{Gigliola Staffilani }

\address{Department of Mathematics, Princeton University}

\email{irod@math.princeton.edu}

\address{Departments of Mathematics, Massachusetts Institute of
Technology and Princeton University}
\curraddr{Department of
Mathematics, Johns Hopkins University}
\email{yanir@member.ams.org}

\address{Department of Mathematics, Massachusetts Institute of
Technology} \email{gigliola@math.mit.edu}

\thanks{\hglue-12pt November 5, 2008. Revised December 2, 2008.}
\thanks{\hglue-12pt Mathematics Subject Classification (2000):
Primary 35Q55.  
Secondary 53C44, 
35B10  
32Q15, 
42B35, 
15A23. 
}
\thanks{\hglue-12pt Keywords: Schr\"odinger flow, periodic NLS, Strichartz estimates, K\"ahler
manifolds.}

\begin{document}

\maketitle

\begin{abstract}

We establish the global well-posedness of the initial value problem for the Schr\"odinger map
flow for maps from the real line into \K manifolds and for maps from the circle into Riemann surfaces.
This partially resolves a conjecture of W.-Y. Ding.

\end{abstract}

\section{Introduction}

In this article we study the \Sch map flow from a one-dimensional
domain into a complete \K manifold. First, we show that when the domain is the real line the flow
exists for all time. Second, we show that when the domain is the circle and the target is a Riemann surface
the flow also exists for all time. The main contribution of this article is to bring Bourgain's work on the periodic
cubic nonlinear \Sch equation (NLS) to bear on the geometric situation at hand.

Let $(M,g)$ be a complete Riemannian manifold of dimension $m$, 
and let
$(N,\o,J,h)$ be a complete symplectic manifold of dimension $2n$ with a compatible almost
complex structure $J$, that is such that $\o(J\cdot,J\cdot)=\o(\cdot,\cdot)$ and such that
$h(\cdot,\cdot)=\o(\cdot,J\cdot)$ defines a complete Riemannian metric on $N$. Associated to this
data is the space of all smooth maps from $M$ to $N$, the Fr\'echet manifold $X:=C^\infty(M,N)$,
endowed with a symplectic structure,
$$
\Omega(V,W)|_u=\int_M\ustar\o(V,W)\dVMg, \quad \forall\; V,W\in T_uX=\Gamma(M,\ustarTN),
$$
where the tangent space to $X$ at a map $u:M\ra N$ is the space of
smooth sections of $\ustarTN\ra M$ and where $\dVMg$ denotes the
volume form on $M$ induced by $g$. The form $\Omega$ is
non-degenerate, i.e., endows $X$ with an injective map $TX\ra
T^\star X$.

Define the energy function on $X$ by
$$
E(u)=\half\int_M |du|_{g^\sharp\otimes\ustarh}^2\dVMg,
$$
where we denote by $g^\sharp$ the metric induced by $g$ on
$T^\star M$ and where we view $du$ as a section of $T^\star
M\otimes\ustarTN\ra M$ and equip this bundle with the metric
$g^\sharp\otimes\ustarh$.

The almost-complex structure on $N$ induces one on $X$ and a corresponding compatible Riemannian metric
$$
G(V,W)|_u=\int_M\ustarh(V,W)\dVMg, \quad \forall\; V,W\in T_uX=\Gamma(M,\ustarTN).
$$
In the infinite-dimensional setting not every function will necessarily have a gradient. However,
if we let $\{u_t\}_{t\in[-1,1]}$ be a smooth
family of maps with $u_0=u$ and denote by $W=\ddto{u_t}$ a variation, then
\begin{equation}
\begin{array}{lll}
\label{EnergyGradientEq}
dE(W)|_{u_0}&=\dis\ddtBigo{E(u_t)}=\int_M g^\sharp\otimes\ustarh(du,d W)\dVMg
\cr\cr
&= -\int_M \ustarh(\tr_{g^\sharp}\N du, W)\dVMg,&
\end{array}
\end{equation}
and hence the gradient of $E$ exists and is given by
\begin{equation}
\label{EnergyGradientDefiningEq}
\N^G E|_u =-\tau(u),
\end{equation}
where
$\tau(u):=\tr_{g^\sharp}\N du$ is called the tension field of $u$ and $\N$ is the connection
on $T^\star M\otimes\ustarTN\ra M$ induced from the Levi-Civita
connection on $(M,g)$ and the pulled-back one from $(N,h)$.
The corresponding gradient flow
\begin{equation}
\label{ESFlowEq}
\dudt=\tau(u),\quad u|_{\{0\}\times M}=u_0,
\end{equation}
is the classical harmonic map flow introduced by Eells and Sampson \cite{ES} that has been
extensively studied.

Now the symplectic gradient of $E$ also exists and is given by
$$
\N^\O E|_u =-J\tau(u).
$$
The corresponding Hamiltonian flow on $(X,\O)$, introduced by Ding-Wang and Terng-Uhlenbeck \cite{DW1,TU},
\begin{equation}
\label{SchEq}
\dudt=-J\tau(u), \quad u|_{\{0\}\times M}=u_0,
\end{equation}
is called the \Sch map flow.

While the energy decreases along (\ref{ESFlowEq}), for
(\ref{SchEq}) the flow is contained in an energy level-set since
for maps of finite energy we have by (\ref{EnergyGradientEq})
\begin{equation}
\label{EnergyConservationEq}
\frac{dE(u(t))}{dt}=
-\int_M \ustarh(\tau(u), \dudt)\dVMg=0.
\end{equation}
For (\ref{ESFlowEq}) one typically
expects to converge to a harmonic representative of the homotopy class of $u_0$ under some
geometric assumptions (e.g., negatively curved target \cite{ES})
while (\ref{SchEq}) seems to be describing some rather very different behavior. Analytically this may be described by the transition from the parabolic equation (\ref{ESFlowEq}) to the borderline case
(\ref{SchEq}) whose symbol has purely imaginary eigenvalues. Note also that for the \Sch flow there
is no preferred time direction.

One problem common to both flows is the question of existence and uniqueness.
Indeed since the flows are defined on infinite-dimensional spaces one cannot expect
global existence/well-posedness in general.%
\footnote{ We note that for the \Sch flow the question of global
existence is equivalent to the existence of a ``symmetry" of
$(X,\Omega)$, i.e. a one-parameter subgroup of Hamiltonian
diffeomorphisms of $(X,\Omega)$ for the energy function $E$
integrating $\N^\Omega E$.} \footnote{ Until recently no results
were known for general symplectic targets; see Chihara \cite{Ch}
for recent work on local well-posedness in this setting.}
Restricting to the \K case, there is a similarity between the two
as far as the local existence is concerned: Results of Ding-Wang
and McGahagan  
show that at least locally (\ref{SchEq}) can be approximated by equations of either parabolic (in the sense
of Petrovski\u\inodot \ (see, e.g., \cite{EZ})) or hyperbolic character. 
As a consequence the following result holds for maps of finite energy.

\begin{theo}
\label{DingWangThm} {\rm(See \cite{DW2,M2}.)} Let $(M^m,g)$ be a
complete Riemannian manifold and let $(N,J,h)$ be a complete \K
manifold with bounded geometry.
For integers $k>m/2+1$ 
equation (\ref{SchEq})
with $u_0\in W^{k,2}(M,N)$
 admits a unique solution $u\in C^0([0,T],W^{k,2}(M,N))$
where $T<T_0$ and $T_0$ depends on $\n{\N u_0}{W^{[m/2]+1,2}}$ and
the geometry of $N$ alone. Moreover, there exist positive
constants $C_1,C_2$ depending only on these quantities such that
$$
\n{\N u(t)}{W^{[m/2]+1,2}}\le C_1/(T_0-t)^{C_2},  \q \all t\in[0,T_0).
$$
In particular, if $u_0\in C^\infty(M,N)$ then $u\in
C^\infty([0,T]\times M,N)$.
\end{theo}

Here by bounded geometry we mean uniform bounds on the injectivity radius and
the curvature tensor and its derivatives. This is
automatically true for compact targets.

The main difficulty lies, therefore, in understanding the global
behavior.

Previous results of a global nature are mostly concerned with the
one-dimensional domain case and are all restricted to the case of
a special target \K manifold. We recall the following
non-exhaustive list of works. The flow on $(S^1,\h{can})\ra (S^2,
\h{can})$ corresponding to the classical model for an isotropic
ferromagnet was studied from the mathematical point of view by
Sulem, Sulem and Bardos \cite{SSB} who obtained local
well-posedness for the initial value problem as well as partial
global results. Zhou, Guo and Tan \cite{ZGT} studied the global
well-posedness problem using a parabolic approximation which was
later put to use by Ding and Wang \cite{DW1} as well as  by Pang,
Wang and Wang \cite{PWW1} to prove global existence and uniqueness
of smooth solutions of maps from $(S^1,\h{can})$ to any constant
sectional curvature \K target as well as to Hermitian locally
symmetric spaces \cite{PWW2} using a conservation law. The latter
also treats the inhomogeneous flow which can be essentially viewed
as the Schr\"odinger flow with domain $S^1$ equipped with a
different metric. Terng and Uhlenbeck \cite{TU} gave a detailed
study of the flow from the Euclidean line into Grassmannians.
Chang, Shatah and Uhlenbeck \cite{CSU} proved existence and
uniqueness of global smooth solutions for maps of the Euclidean
line into a compact Riemann surface. In addition, they treated
maps of the Euclidean plane into a compact Riemann surface under
the assumption of small initial energy and certain symmetries.
Finally, see Bejenaru, Ionescu, Kenig and Tataru \cite{BIK, BIKT}
for recent work on global well-posedness in the case of maps from
Euclidean space into $(S^2,\h{can})$ under a certain smallness
assumption.

Note that in all of these results one restricts the target to a rather small class of \K manifolds.
We recall the following conjecture of Ding.

\begin{conj} {\rm(See \cite{D}.)}
\label{DingConjecture}
The \Sch map flow is globally well-posed for maps from one-dimensional domains into compact \K manifolds.
\end{conj}

The main results of this article are a partial answer to this conjecture. Namely, we
establish the global well-posedness of the one-dimensional
\Sch flow into general \K manifolds when the domain is the real line, and into Riemann surfaces when
the domain is the circle.

\begin{theo}
\label{MainThmOne}
Let $(M,g)=(\RR,dx\otimes dx)$, let $(N,J,h)$ be
a complete \K manifold with bounded geometry, 
and let $k\ge 2$ be an integer. 
The flow equation (\ref{SchEq})
with $u_0\in W^{k,2}(\R,N)$ admits a unique solution $u\in C^0(\RR,W^{k,2}(\RR,N))$.
In particular $u$ is smooth if $u_0$ is. 
\end{theo}

\begin{theo}
\label{MainThmTwo}
Let $(M,g)=(S^1,dx\otimes dx)$, let $(N,J,h)$ be
a Riemann surface with bounded geometry, 
and let $k\ge 2$ be an integer. 
The flow equation (\ref{SchEq})
with $u_0\in W^{k,2}(S^1,N)$ admits a unique solution $u\in C^0(\RR,W^{k,2}(S^1,N))$.
In particular $u$ is smooth if $u_0$ is. 
\end{theo}

\begin{remark} {\rm
From the physical point of view, the \Sch map flow may also be introduced
as a generalization of the Heisenberg model for a ferromagnetic spin system. The classical
model for this physical system precisely corresponds to maps from the standard
circle into
$N=S^2$ with the standard metric and complex structure \cite{LL} (for some background see,
e.g., \cite{D,DW1,M1,SSB}).
Perhaps the most physically natural generalization of the classical model would be to vary
the metric on the target $S^2$, however it seems that even for small perturbations of the
round metric on $S^2$ global well-posedness was not known before.
Theorem \ref{MainThmTwo} establishes
the global well-posedness of the Cauchy problem describing this physical model
when the metric on $S^2$ is arbitrary.
}
\end{remark}

Let us now outline the key ideas of the proofs.
According to Theorem \ref{DingWangThm} we have existence of a local (in time) solution.
The strategy of the proof is: First, we translate the flow equation
into a system of nonlinear \Sch (NLS) equations.
Second, for {\it this} system of equations we obtain
an a priori estimate (and hence local well-posedness)
in a weaker norm than that in Theorem \ref{DingWangThm}, namely in an
appropriate Strichartz norm for $M=\RR$ and in $L^4$ for $M=S^1$.
These estimates
are crucial since they only depend on the initial energy (that is a conserved quantity)
and that in a manner that can be readily converted into a global a priori estimate (and hence
global well-posedness for the system of NLS equations) in the same
space. Taking derivatives of the flow equation we obtain global a priori estimates
in stronger norms and these in turn may be converted back to imply global well-posedness
for our original Cauchy problem in $W^{k,2}$ for all $k\ge2$.
While the proofs of both Theorem \ref{MainThmOne} and Theorem \ref{MainThmTwo}
follow the same general scheme, nevertheless  there are substantial differences between the two as we now explain.

First, we consider the case of the real line. This case is
considerably simpler due to the non-compactness and
dispersiveness. Here we follow Chang-Shatah-Uhlenbeck and write
the flow equation in terms of a parallel frame observing the \K
condition allows one to readily generalize their computations from
the Riemann surface case to higher dimensions. The flow equation
then reduces to a system of NLS equations. The same Strichartz
type calculations as in their study of the Riemann surface case
then apply.

Second, we treat the case $M=S^1$. This case is considerably more
difficult and is the main contribution of this article. There are
two main difficulties. First, using a parallel frame introduces
holonomy and so the resulting NLS equations live on $\RR$, the
universal cover of $S^1$, instead of on $S^1$ itself. To overcome
this we use a certain space-time transformation in order to obtain
a system of NLS equations on $S^1$ in terms of the holonomy
representation of $N$. In addition we need to estimate the
variation of the holonomy along the flow. Second, since $S^1$ is
compact the equations are no longer dispersive. To overcome that
we adapt Bourgain's results on the cubic NLS to our setting in
order to prove local well-posedness in $L^4$ that depends in such
a way on the initial data that it may be used to obtain global
well-posedness in the same space.

The article is organized as follows. In Section \ref{SectionR} we
treat the case of the real line. In Section
\ref{SectionCircleOneD} we treat the case of maps from the circle
into a Riemann surface. Finally, in Section
\ref{SectionCircleHigherD} we discuss some of the difficulties
that arise when trying to apply our approach to treat maps from
the circle into higher dimensional \K manifolds. It is conceivable
that some of these ideas might be related to showing finite time
blow-up for higher-dimensional domains.

The authors thank B. Dai for a careful reading and for correcting
an error in the computation of the holonomy in an earlier version
of this article. This material is based upon work supported in
part under a National Science Foundation Graduate Research
Fellowship and grants DMS-0406627,0702270,0602678. Y.A.R. was also
partially supported by a Clay Mathematics Institute Liftoff
Fellowship.

\section{Maps of the real line into a \K manifold}
\label{SectionR}

In this section we consider maps from the Euclidean real line
into a complete \K manifold $(N,J,h)$ of complex dimension $n$.

The \Sch equation (\ref{SchEq}) becomes

\begin{equation}
\label{SchREq}
J\Nt u-\Nx\Nx u =0,\quad u(0)=u_0,
\end{equation}
\begin{equation}
\label{FiniteEnergyEq}
E(u_0) =\half\int_\R |du_0|_{\ppx\otimes\ppx\otimes\ustarh}^2dx<\infty,
\end{equation}
where we use the abbreviated notations
$\Nt=\nabla_{u_{\star}\ppt}, \Nx=\nabla_{u_{\star}\ppx}$ and $\Nt
u=u_{\star}\ppt=\frac{\pa u}{\pa t}, \Nx u=u_{\star}\ppx=\frac{\pa
u}{\pa x}$ and denote derivatives of a function $f$ by $f_{,x}$
and $f_{,t}$. The key idea in this section, going back to Chang,
Shatah and Uhlenbeck \cite{CSU}, is to rewrite (\ref{SchREq}) in
an appropriate frame along the image in such a way that
(\ref{SchREq}) reduces to a system of nonlinear \Sch (NLS)
equations. In fact our proof closely follows their approach for
the Riemann Surface case observing that it readily generalizes to
\K targets of arbitrary dimension.

Now assume that $u:I\times \R\ra N$ is a solution of
(\ref{SchREq}) with $I$ a neighborhood of $0$ in $\R$ (given, for
example, by Theorem \ref{DingWangThm}). Choose an orthonormal
frame $\{e_1,\ldots,e_{2n}\}$ for $u^\star TN$ with respect to
$h$. We further reduce the structure group to $U(n)\subseteq
O(2n)$ by letting $e_{n+1}=Je_1,\ldots,e_{2n}=Je_n$. We identify
$U(n)$ with its image in $O(2n)$ under the map $\iota:GL(n,\C)\ra
GL(2n,\R)$ given by
$$
\iota(A+\i B)
=
\begin{pmatrix}
A&-B\cr B&A
\end{pmatrix}.
$$
Note
that if
$v=x+\i y\in\C^n$ and
$$
\iota(v)
=
\begin{pmatrix}
x\cr y
\end{pmatrix}
$$
then
$$
\iota(Av)=\iota(A)\iota(v).
$$
We will use this identification frequently, sometimes omitting the
reference to the map $\iota$. In the following we let latin and
greek indices take values in $\{1,\ldots,2n\}$ and
$\{1,\ldots,n\}$, respectively. For both alphabets we use the
notation
$$
\bar\cdot=\cdot+n-1 \;(\h{mod\ }
2n)\;+1.
$$
Therefore barred greek indices take values in $\{n+1,\ldots,2n\}$.
Put $e_{\b j}=Je_j$ so $e_{\b\a}=Je_\a,\; e_{\overline{\a+n}}=-e_\a$.
We let, for example, $L^p(\RR_x)$ and $L^p(\RR_t)$ denote the spaces $L^p(\RR,dx)$ and $L^p(\RR,dt)$,
respectively. Finally, given a map $u:(M,g)\ra (N,h)$ by $u\in W^{k,p}(M,N)$
we will mean that $\sum_{j=0}^{k-1}\big|\big||\N^j du|\big|\big|_{L^p}<\infty$.
For example, in this notation we have $E(u)=||u||^2_{W^{1,2}(M,N)}$.

Now we may view the flow equation (\ref{SchREq}) in this frame.
Write for each $(t,x)\in I\times\R$
\begin{equation}
\label{VectorFieldsInTermsOfBasisEq}
\Nx u=\sum_{j=1}^{2n}h(\Nx u,e_j)e_j=:a^je_j,\qquad
\Nt u=\sum_{j=1}^{2n}h(\Nt u,e_j)e_j=:b^je_j,
\end{equation}
where we make use of the Einstein summation convention, namely the appearance of an index both as a subscript
and a superscript indicates summation.

Equation (\ref{SchREq}) can be rewritten as
\begin{equation}
\label{FlowEqInCoordinatesEq}
b^je_{\b j}-a^j_{,x}e_j-a^j\Nx e_j=0.
\end{equation}
The conservation of energy (see (\ref{EnergyConservationEq})) is expressed as
\begin{equation}
\label{EnergyConservationSecondEq}
E(u(t))=E(u_0)=\half\int_\R\sum_{l=1}^n (a_l)^2 dx=||a(t)||_{L^2(\RR_x)}^2,\qq\all\,t\in\R.
\end{equation}
Note that
\begin{equation}
\label{BracketVectorFieldsEq}
0=u_\star\Big[\ppt,\ppx\Big]=[\N_t u,\N_x u].
\end{equation}
Since $\N J=0$, differentiating (\ref{SchREq}) in space yields
$$
J\Nx\Nt u-\Nx\Nx\Nx u=0,
$$
which becomes, using
(\ref{BracketVectorFieldsEq}) and (\ref{VectorFieldsInTermsOfBasisEq}),
\begin{equation}
J(a^j_{,t}e_j+a^j\Nt e_j)-a^j_{,xx}e_j-2a^j_{,x}\Nx e_j-a^j\Nx\Nx e_j=0.
\end{equation}

We impose the gauge-fixing condition
\begin{equation}
\label{GaugeEq}
\N_x e_j=0,\quad j=1,\ldots,2n.
\end{equation}
The resulting frame along the image is still unitary since the
complex structure commutes with parallel transport. Equation
(\ref{FlowEqInCoordinatesEq}) becomes
\begin{equation}
\label{BjInTermsOfAj}
b^j=a^{\b j}_{,x}.
\end{equation}
Note that $u^\star TN\ra \R$ is trivial and that (\ref{GaugeEq})
amounts to fixing a trivializing parallel frame. With this choice
the flow on $\ustarTN$ is given by
\begin{equation}
\label{FlowAfterGaugeEq}
a^j_{,t}e_{\b j}-a^j_{,xx}e _j=-a^j\Nt e_{\b j}.
\end{equation}

Along the image, using (\ref{GaugeEq}) and (\ref{BjInTermsOfAj}), and letting $R$ denote the curvature tensor
of $(N,h)$, we have
\begin{equation}
\begin{array}{lll}
\dis\N_x\N_te_{\b j} &= R(\N_xu,\N_tu)e_{\b j}= a^kb^lR_{kl{\b
j}}^{\phantom{kl\b j}q}e_q \cr\cr \dis &=(a^\a b^\be R_{\a\be{\b
j}}^{\phantom{\a\be\b j}q}+ a^{\b\a}b^\be R_{\b\a \be{\b
j}}^{\phantom{\a\be\b j}q}+ a^{\a}b^{\b\be} R_{\a \b\be{\b
j}}^{\phantom{\a\be\b j}q}+ a^{\b\a}b^{\b\be} R_{\b\a \b\be{\b
j}}^{\phantom{\a\be\b j}q} )e_q\cr\cr \dis &= (a^\a a^{\b\be}_{,x}
R_{\a\be{\b j}}^{\phantom{\a\be\b j}q}+ a^{\b\a}a^{\b\be}_{,x}
R_{\b\a \be{\b j}}^{\phantom{\a\be\b j}q}- a^{\a}a^{\be}_{,x}
R_{\a \b\be{\b j}}^{\phantom{\a\be\b j}q}- a^{\b\a}a^{\be}_{,x}
R_{\b\a \b\be{\b j}}^{\phantom{\a\be\b j}q} )e_q\cr\cr \dis &=
\dis\sum_{\a,\be}\big[(a^\a a^{\b\be})_{,x} R_{\a\be{\b
j}}^{\phantom{\a\be\b j}q}+
\half[(a^{\b\a}a^{\b\be})_{,x}+(a^{\a}a^{\be})_{,x}] R_{\b\a
\be{\b j}}^{\phantom{\a\be\b j}q}
\big]e_q,\cr\cr
\end{array}
\end{equation}
where we have used the \K condition once more:
$$
R_{\a\be\b j}^{\phantom{\a\be\b j}q}
=
R_{\b\a\b\be\b j}^{\phantom{\a\be\b j}q},\;
R_{\b\a\be\b j}^{\phantom{\a\be\b j}q}
=
-R_{\a\b\be\b j}^{\phantom{\a\be\b j}q}.
$$
Equation (\ref{FiniteEnergyEq}) implies that $\lim_{x\ra\pm\infty}a^i(t,x)=0$.
Therefore, since
$$
\N_x h(\N_te_{\b j},e_q)=h(\N_x\N_te_{\b j},e_q),
$$
we have
\begin{equation}
\begin{array}{lll}
&\dis\N_te_{\b j}(t,x) = \cr\cr
&\dis\q\;\sum_{q=1}^{2n}h(\N_te_{\b
j},e_q)(t,-\infty)e_q(t,x)\cr\cr
&\q\;\dis+\Big[\sum_{\a,\be}\big[a^\a a^{\b\be} R_{\a\be{\b
j}}^{\phantom{\a\be\b j}q}+ \half[a^{\b\a}a^{\b\be}+a^{\a}a^{\be}]
R_{\b\a \be{\b j}}^{\phantom{\a\be\b j}q} \big](t,x)\cr\cr
&\q\;\dis-\int_{(-\infty,x]}\sum_{\a,\be}\big[a^\a a^{\b\be}
R_{\a\be{\b j}}^{\phantom{\a\be\b j}q}{}_{,x}+
\half(a^{\b\a}a^{\b\be}+a^{\a}a^{\be}) R_{\b\a \be{\b
j}}^{\phantom{\a\be\b j}q}{}_{,x}
\big](t,y)dy\Big]e_q(t,x) \cr\cr
& \qquad\qquad\;\;\!=: A_{\b j}^q(t,-\infty) e_q(t,x)+[P_{\b
j}^{q}(t,x)+Q_{\b j}^q(t,x)]e_q(t,x).\cr
\end{array}
\end{equation}
Note that by using (\ref{FlowEqInCoordinatesEq}) we have
\begin{equation}
\label{ChristoffelTermEq}
\N_t e_{\b j}=b^k\Gamma^p_{k\b j} e_p=a^{\b k}_{,x}\Gamma^p_{k\b j} e_p.
\end{equation}
Hence,
$h(\N_te_{\b j},e_q)=a^{\b k}_{,x}\Gamma^q_{k\b j}$ and so we may assume that $A_{\b j}^q$ vanishes at
$(t,-\infty)$. To justify this note that this is indeed the case for the local solution
of our equation given by Theorem \ref{DingWangThm}; even though
this assumption makes use of the finiteness of the $W^{2,2}$ norm
of that local solution, the important point is that eventually our estimates will not depend
on the $W^{2,2}$ norm of $u$ (equivalently on the $W^{1,2}$ norm of $a$),
and so the proof of the a priori estimate for the system of NLS equations
(\ref{NLSSystemREq}) below (for $a$) goes through, with this assumption.

\noindent
Note that $R_{klp}^{\phantom{klp}q}{}_{,x}=a^sR_{klp}^{\phantom{klp}q}{}_{,s}$.
Therefore
$|P_{j}^{q}|< C||\bfa||^2$ and $|Q_{ j}^{q}|< C\int_M||\bfa||^3dx$.
Finally
(\ref{FlowAfterGaugeEq}) transforms to the following system of NLS equations
\begin{equation}
-a^{\b\ga}_{,t}-a^\ga_{,xx}=-a^jP_{ j}^\ga-a^jQ_{ j}^\ga, \qquad \gamma=1,\ldots,n,
\end{equation}
\begin{equation}
a^{\ga}_{,t}-a^{\b\ga}_{,xx}=-a^jP_{ j}^{\b\ga}-a^jQ_{ j}^{\b\ga}, \qquad \gamma=1,\ldots,n,
\end{equation}
or, letting $J_0=\iota(\i I)$,
\begin{equation}
\label{NLSSystemREq}
J_0\bfa_{,t}=\bfa_{,xx}-\bfP\cdot\bfa-\bfQ\cdot\bfa.
\end{equation}
where $\bfa=(a^1,\ldots, a^{2n})^T,\; \bfP=(P_j^k), \; \bfQ=(Q_j^k)$.

Equivalently, using the aformentioned identification $\iota$ of $GL(n,\C)$ with a subset of $GL(2n,\R)$,

\begin{equation}
\label{NLSEq}
\i\bfPhi_{,t}=
\bfPhi_{,xx}-
\bfS\cdot\bfPhi-\bfT\cdot\bfPhi,
\end{equation}

\smallskip
\noindent
where
$
\bfPhi=\iota^{-1}(\bfa)=(a^1+\i a^{\b 1},\ldots,a^n+\i a^{\b n})^T, \;
\bfS=(S_\a^\be)=\iota^{-1}(\bfP),\;
\bfT=(T_\a^\be)=\iota^{-1}(\bfQ)
$
and
$|S_{\a}^{\be}|< C||\bfPhi||^2$ and $|T_{\a}^{\be}|< C\int_M||\bfPhi||^3dx$. Here $C$ depends only
on
the geometry of $(N,J,h)$, which we assume to be bounded.

\begin{remark}
\label{RemarkGeneralMetric} {\rm
In the case of a variable
complete smooth metric $(M,g)=(\R,\a^{-1}dx\otimes dx)$ with
$\a>0$ equation (\ref{SchEq}) becomes
\begin{equation}
b^ke_{\b k}=\a a^k_{,x}e_k+\half\a_{,x}a^ke_k,
\end{equation}
which can then be transformed, as before, to
\begin{equation}
J_0\bfa_{,t}
=
\a\bfa_{,xx}+\frac{3\a_{,x}}2\bfa_{,x}+\frac{\a_{,xx}}2\bfa
-\bfP\cdot\bfa-\bfQ\cdot\bfa.
\end{equation}
Equivalently, again using the map $\iota:GL(n,\C)\ra GL(2n,\R)$,
\begin{equation}
\label{NLSVariableMetricEq}
\i\bfPhi_{,t}=
\a\bfPhi_{,xx}+
\frac{3\a_{,x}}2\bfPhi_{,x}+
\frac{\a_{,xx}}2\bfPhi-
\bfS\cdot\bfPhi-\bfT\cdot\bfPhi.
\end{equation}
The only obstacle to treating this equation using the methods
below is the first derivative term on the right hand side.}
\end{remark}
\medskip

Therefore we have reduced the original flow equation for the map to a system of NLS equations for
the frame coefficients of the gradient of the map.

We are now in a position to demonstrate the following theorem that is the main result of this section.

\begin{theo}
\label{RThm}
Let $(M,g)=(\R,dx\otimes dx)$ and let $(N,J,h)$ be
a complete \K manifold with bounded geometry. 
Then for integers $k\ge 2$ 
equation (\ref{SchEq})
with $u_0\in W^{k,2}(\R,N)$ admits a unique solution $u\in C^0(\R,W^{k,2}(\R,N))$.
In particular $u$ is smooth if $u_0$ is.
\end{theo}

\begin{proof}
First, we have by Theorem \ref{DingWangThm} local well-posedness
of the original Schr\"odinger map flow (\ref{SchREq}) in
$C^0(\RR_t,W^{k,2}(\RR,N))$ for $k\ge2$. Our purpose is now to
prove a local a priori estimate in a space that is
weaker than $W^{2,2}(\RR,N)$. Equivalently, we will prove an
estimate for the frame coefficients $a$ in a norm weaker than $W^{1,2}(\RR,\RR^{2n})$.
More precisely, we first assume we are
given a smooth initial data and prove estimates on the smooth
local solution given by Theorem \ref{DingWangThm} in the weaker
space $L^4(\RR_{t,loc},W^{1,\infty}(\RR,N))\cap
C^0(\RR_t,W^{1,2}(\RR,N))$ that will be independent of this
smoothness assumption. The key will be that these estimates will
depend only on the initial energy $E(u_0)$. Since the energy is
conserved these will then lead to a global a priori estimate in the
same space. Then, by taking derivatives of the flow equation we
will obtain the global well-posedness claimed in the statement.

Each of the equations in the system (\ref{NLSEq}) is of exactly the same type as the equation
obtained by Chang-Shatah-Uhlenbeck \cite{CSU} for the case of a target Riemann surface. The only difference
is that each $\bfPhi^j$ depends also on the other $\bfPhi^k, k=1,\ldots,n$. However this dependence is only
in the nonlinear terms and not in the terms involving derivatives.
It therefore follows that their proof of local well-posedness in the Strichartz spaces
$L^4(\RR_t,L^\infty(\RR_x))\cap C^0(\RR_t,L^2(\RR_x))$ (see (\ref{LpLqNormNotationEq}) below
for notation)
carries over to our system of NLS equations.
Since this estimate depends only on the initial energy (that is a conserved quantity) we obtain
global well-posedness in $L^4(\RR_t,L^\infty(\RR_x))$. Taking a derivative of the equation
and working in the intersection of the Strichartz spaces $L^4(\RR_t,L^{1,\infty}(\RR_x))$
and $L^\infty(\RR_t, W^{1,2}(\RR_x))$ as in their original proof then yields global
well-posedness in $W^{2,2}(\RR,N)$ for the original flow equation
(i.e., global existence in time for initial data in $W^{2,2}(\RR,N)$).
One may then show using similar computations global well-posedness in $W^{k,2}$ for
each $k\ge 2$.

In fact, although we will not carry this out here, one may prove local (and hence global) well-posedness
for (\ref{NLSEq}) in other Strichartz spaces (these are by definition the spaces
$L^q(\RR_t,L^r(\RR_x))$ specified by Lemma \ref{StrichartzLemma} below) as well, e.g., $L^6(\RR\times\RR)$.

For the benefit of the reader that may not be familiar with standard Strichartz estimate techniques
we include here the detailed proof of the Chang-Shatah-Uhlenbeck  $L^4(\RR_t,L^\infty(\RR_x))$ estimate. No originality
is claimed here.
This also serves to provide some perspective on the differences between this
case and the case of the circle, treated in the following sections.
The estimates on higher derivatives are similar and for that we refer to their original article.

Suppose a function $c:\R\times\R\ra\C$ satisfies the NLS equation
\begin{equation}
\label{NLSsingleEq}
\i\, c_{,t}=c_{,xx}+F, \q \all t\in[0,T],\qq c(0)=f,
\end{equation}
for some function $F:[0,T]\times\R\ra\C$ that
may depend on $c$ nonlinearly (but not on its derivatives).
One then has the integral expression (Duhamel formula)
\begin{equation}
\label{DuhamelEq}
c(t,x)=
\int_{\R}f(y){e^{-\i|x-y|^2/4t}\ov\sqrt{2\pi t}} dy
-\i\int_0^t \int_{\R} F(s,y){e^{-\i|x-y|^2/4(t-s)}\ov\sqrt{2\pi (t-s)}} dy\w ds.
\end{equation}

Denote the \Sch operator by
$$
S(t)f:=\int_{\R}f(y){e^{-\i|x-y|^2/4t}\ov\sqrt{2\pi t}} dy.
$$

We now recall the Strichartz estimates (on $\RR$). For appropriate $q,r$ we denote by $L^q(\RR,L^r(\RR))$
the Banach space equipped with the norm
\begin{equation}
\label{LpLqNormNotationEq}
||f||_{L^q(\RR,L^r(\RR))}
:=
\big|\big| ||f||_{L^r(\RR_x)}\big|\big|_{L^q(\RR_t)}.
\end{equation}

\begin{lem} {\rm (See \cite{C}, p. 33.)}
\label{StrichartzLemma}
Let $q,r$ satisfy $\frac2q+\frac1r=\frac 12$, with $r\in[2,\infty]$ and let $f\in L^2(\RR_x)$.
Then the function $t\mapsto S(t)f$ belongs to $L^q(\RR_t,L^r(\RR_x))\cap C^0(\RR_t,L^2(\RR_x))$
and there is a constant independent of $(q,r)$ and $f\in L^2(\RR)$ such that
\begin{equation}
\label{StrichartzEstimateEq}
||S(\cdot)f||_{L^q(\RR_t,L^r(\RR_x))}\le C||f||_{L^2(\RR)}.
\end{equation}
\end{lem}

In our situation we know that $||\bfPhi||_{L^2}$ is constant in time (recall (\ref{EnergyConservationSecondEq})).
Assume also that
$\bfPhi\in L^4([0,T],L^\infty(\RR_x))$. We will now show that in fact the $L^4([0,T],L^\infty(\RR_x))$ norm
of $\bfPhi$ is controlled by its $L^2$ norm and the geometry.  This will imply local and eventually
global well-posedness in $L^4(\RR_t,L^\infty(\RR_x))$.

Let $F=-\bfS\cdot\bfPhi-\bfT\cdot\bfPhi$. Then
\begin{equation}
\label{PhijEvolutionEq}
\bfPhi^j(t)
=
S(t-t_1)\bfPhi^j(t_1)
-
\i \int_{t_1}^t S(t-s)F(s,\cdot) ds.
\end{equation}
The first term of (\ref{PhijEvolutionEq}) is in
$L^4([t_1,t],L^\infty(\RR))$ by the Strichartz estimate
(\ref{StrichartzEstimateEq}). We will now show that the second
term is also in this space.

First, we consider the term $\bfS\cdot\bfPhi\le C||\bfPhi||^3$.
By applying the Strichartz estimate under the integral sign and using energy conservation we obtain
\begin{equation}
\label{AnotherInequalityNonlinearTermEq}
\begin{array}{lll}
&\dis\big|\big|\int_{t_1}^t S(t-s)(\bfS\cdot\bfPhi)(s,\cdot) ds\big|\big|_{L^4([t_1,t_2],L^\infty(\RR))}
\cr\cr
&\qquad\qquad\le
\dis\int_{t_1}^t \big|\big| S(t-s)(\bfS\cdot\bfPhi)(s,\cdot)\big|\big|_{L^4([t_1,t_2],L^\infty(\RR))} ds
\cr\cr
&\qquad\qquad\le
C \dis\int_{t_1}^t \big|\big| (\bfS\cdot\bfPhi)(s,\cdot)\big|\big|_{L^2(\RR)} ds
\cr\cr
&\qquad\qquad\le
C' \dis\int_{t_1}^t \big|\big| |\bfPhi(s,\cdot)|^3\big|\big|_{L^2(\RR)} ds
\cr\cr
&\qquad\qquad\le
C' \dis\int_{t_1}^t \big|\big| \bfPhi(s,\cdot)\big|\big|^2_{L^\infty(\RR)}
\big|\big| \bfPhi(s,\cdot)\big|\big|_{L^2(\RR)} ds
\cr\cr
&\qquad\qquad\le
C'' \dis\int_{t_1}^t \big|\big| \bfPhi(s,\cdot)\big|\big|^2_{L^\infty(\RR)}ds
\le
C'' |t-t_1|^{1/2} \big|\big|\bfPhi\big|\big|^2_{L^4([t_1,t_2],L^\infty(\RR))}.
\end{array}
\end{equation}

Second, we consider the term $\bfT\cdot\bfPhi\le
C||\bfPhi||\int_{\RR}||\bfPhi||^3dx$. Again, by applying the
Strichartz estimate under the integral sign and using energy
conservation we obtain

\begin{equation}
\begin{array}{lll}
\label{InequalitySecondNonlinearTermEq}
&\dis\big|\big|\int_{t_1}^t S(t-s)(\bfT\cdot\bfPhi)(s,\cdot) ds
\big|\big|_{L^4([t_1,t],L^\infty(\RR))}
\cr\cr & \dis\qquad\qquad
\le C' \int_{t_1}^{t} \Big|\Big|\bfPhi
\big|\big||\bfPhi|^3\big|\big|_{L^1(\RR)}\Big|\Big|_{L^2(\RR)}ds.
\cr\cr & \dis\qquad\qquad \le C' \int_{t_1}^{t} \Big|\Big|\bfPhi
\big|\big|\bfPhi\big|\big|_{L^\infty(\RR)}
\big|\big|\bfPhi\big|\big|^2_{L^2(\RR)}\Big|\Big|_{L^2(\RR)}ds.
\cr\cr & \dis\qquad\qquad \le C'' \int_{t_1}^{t}
\big|\big|\bfPhi\big|\big|_{L^\infty(\RR)}ds \le
C''|t-t_1|^{3/4}\big|\big|\bfPhi\big|\big|_{L^4([t_1,t],L^\infty(\RR))}.
\end{array}
\end{equation}

\bigskip

Combining (\ref{AnotherInequalityNonlinearTermEq}) and (\ref{InequalitySecondNonlinearTermEq})
we thus obtain, by choosing $|t-t_1|$ small enough (depending only on the initial energy), an estimate
on $||\bfPhi||_{L^4([t_1,t],L^\infty(\RR))}$, depending only on the geometry of $(N,h)$ and the initial
energy. This then implies global well-posedness in $L^4(\RR_t,L^\infty(\RR_x))$,  
by which we mean that for all $T>0$ we have
$$
\int_0^T ||\bfPhi(s,\cdot)||^4_{L^\infty(\RR_x)}ds<\infty.
$$
As indicated earlier the higher derivatives estimates follow similar computations. This concludes
the proof of Theorem \ref{RThm}.
\end{proof}

\begin{remark}
{\rm It is essential to use Theorem \ref{DingWangThm} since even
after we reduce the \Sch map flow to a system of NLS equations and
after proving that a unique global solution for (\ref{NLSEq})
exists it is not completely obvious how to go from such a solution
for $\N_x u$ to an actual map $u$ into $N$. }
\end{remark}

\section{Maps of the circle into a Riemann surface}
\label{SectionCircleOneD}

In this section we consider the \Sch map flow with the domain
being the round circle. Compared with the previous section, the
discussion here is more delicate due to the fact that the domain
is no longer simply-connected (introduces holonomy) nor
non-compact (lack of dispersion).

Let
$u:I\times S^1\ra N$ where $I\subset\R$ is a neighborhood of $0$.
The bundle
$\ustarTN\ra I\times S^1$ is no longer trivial and so
fixing a frame satisfying (\ref{GaugeEq}) does not yield a trivialization. To describe the
solution of
(\ref{GaugeEq}) 
we work instead with $\Z$-invariant objects over $\R$. We
therefore make the identifications
\begin{equation}
\label{ZZinvariantEq}
\h{Maps}(S^1,N)\isom\h{Maps}(\R,N)^\Z,\quad
\Gamma(I\times S^1,\ustarTN)\isom \Gamma(I\times \R,\ustarTN)^\Z,
\end{equation}
the superscript denoting $\Z$-invariant objects, and
take the freedom to use either one of these identifications interchangibly.
 Similar identifications will be made for
all the other tensor bundles encountered over $I\times S^1$
(e.g., $\ustar (\TstN\otimes\TstN\otimes\TstN\otimes TN)$).

Recall that parallel transport is defined
as a map $P:u\mapsto\h{Aut}(u(0)^\star TN,u(1)^\star TN)$
for all $u\in C^\infty([0,1],N)$, which
on any \K manifold restricts to an operator $P:C^\infty((S^1,\h{pt}),N)\ra \iota(U(n))$
on base-pointed loops. Formally, a solution of (\ref{GaugeEq}) is given by
$e(t,x)=P(u(t)|_{[0,x]})e(t,0)$.
This can be described somewhat more explicitly as follows.

Let
$U$ denote a contractible open set in $N$ and
and let $e_1,\ldots,e_n,e_{n+1}=Je_1, \ldots,$ $e_{2n}=Je_n$ denote a local orthonormal frame.
Assume
$u:I\times\R\ra N$ is a solution of (\ref{SchEq}), a collection of
loops in $N$
which we will initially
assume to be contained
in $U$ (and so, in effect, these loops are all contractible in $N$).
Along the image of our flow we denote by $\a^1,\ldots,\a^{2n}$
the dual 1-forms to $e_1,\ldots,e_{2n}$.
The Levi-Civita
connection along our flow restricted to this patch is represented by a section $A_U=\T kij \alpha^i$
of $\TstarN|_U\otimes u(n)$ which pulls back to a connection form
$u^\star A_U=\T kij a^idx=:B_U dx$
for the pulled-back bundle. A section
$e=E^j e_j$ of the pulled-back bundle (as in (\ref{ZZinvariantEq})) is then
(locally) parallel when
\begin{equation}
\label{GaugeMatrixOneFormEq}
0=\N e=\parr{E^j}x e_j\otimes dx+B_U\cdot e\otimes dx
      =(E^j_{,x}+{B_U}^j_kE^k)e_j\otimes dx.
\end{equation}

The solution of this first-order matrix equation simplifies
considerably in the case $n=1$.
The matrices $B_U$ then lie in the
trivial Lie algebra $so(2)\isom u(1)$ and so their exponentials
commute. One may therefore integrate (\ref{GaugeMatrixOneFormEq})
to obtain
\begin{equation}
\label{FirstHolonomyReprEq}
e(t,1)=\exp(-\int_0^1 B_U dx)e(t,0).
\end{equation}
If $D_u$ is the
disc bounded by $u$ and contained in $U$, and $K$ denotes the
Gaussian curvature of $N$, then Stokes' Theorem gives
$$
e(t,1)=\exp(-\int_{D_u}
dA_U)e(t,0)=\exp(-\int_{D_u}K\dVNh) e(t,0)
$$
(possibly up to a factor of $2\pi$, depending on conventions) from
which it becomes evident that one may relax the assumption above
(for the moment still restricting to contractible loops) and work
globally (one might have two choices for $D_u$ then). Also, we see
that the holonomy factor is independent of the starting point on
the loop. In fact this last fact is seen to be true also for
non-contractible loops. We have therefore a well-defined holonomy
map
$$
P:C^\infty((S^1,\pt),N)\ra
SO(2)=\iota(U(1)).
$$

Next, for general $u$, since $u(0,S^1)$ and $u(t,S^1)$ are homotopic for any
$t\in I$
we may define the surface $D_u=u([0,t]\times S^1)$ and as chains on $N$ $\del D_u=u(t,S^1)-u(0,S^1)$.
Let $K$ denote the Gaussian curvature of $(N,h)$. Then we
have once again by Stokes' Theorem
$$
e(t,1)=P(u)e(t,0)=\exp(-\int_{D_u}K\dVNh)P(u_0) e(t,0),
$$
or for any $x\in\R$ and $l\in\NN$
\begin{equation}
\label{HolonomyEq}
e(t,x+l)=P(u)^le(t,x)=\exp(-l\int_{D_u}K\dVNh)P(u_0)^l e(t,0).
\end{equation}

Therefore a solution of (\ref{GaugeEq}) produces a parallel section of
$\Gamma(\R,\ustarTN)$ rather than of
$\Gamma(\R,\ustarTN)^\Z$. In expressing our $\Z$-invariant tensors in terms of the
frame $\{e_j\}_{j=1}^{2n}$ we therefore use coefficients satisfying a relation appropriately
proportional to (\ref{HolonomyEq}). For example if $v\in\Gamma(\R,\ustarTN)^\Z$ then
we may write $v=v^je_j$ with $v^j(x+l)=P(u)^{-l}v^j(x)$ (while on the other hand sections
of endomorphism tensor bundles require no adjustment when $n=1$).

Going through the computations of \S2 it follows that
Equation (\ref{FlowAfterGaugeEq}) still holds.
We then obtain
\begin{equation}
\begin{array}{lll}
&\mskip-15mu\N_te_{\b j}(t,x)
=\cr\cr
&\q\qq\;\dis
\sum_{q=1}^{2n}h(\N_te_{\b j},e_q)(t,x_0)e_q(t,x)
\cr\cr
&\q\qq\;\dis+\bigg[\sum_{\a,\be}\Big(\big[a^\a a^{\b\be} R_{\a\be{\b j}}^{\phantom{\a\be\b j}q}+
\half[a^{\b\a}a^{\b\be}+a^{\a}a^{\be}] R_{\b\a \be{\b j}}^{\phantom{\a\be\b j}q}
\big](t,x)
\cr\cr
&\q\qq\;\qq\q\dis-
\big[a^\a a^{\b\be} R_{\a\be{\b j}}^{\phantom{\a\be\b j}q}+
\half(a^{\b\a}a^{\b\be}+a^{\a}a^{\be}) R_{\b\a \be{\b j}}^{\phantom{\a\be\b j}q}
\big](t,x_0)\Big)
\cr\cr
&\q\q\qq\q\dis-\int_{[x_0,x]}\sum_{\a,\be}\big[a^\a a^{\b\be}
R_{\a\be{\b j}}^{\phantom{\a\be\b j}q}{}_{,x}+
\half(a^{\b\a}a^{\b\be}+a^{\a}a^{\be}) R_{\b\a \be{\b j}}^{\phantom{\a\be\b j}q}{}_{,x}
\big](t,y)dy\bigg]e_q(t,x)
\end{array}
\end{equation}
The terms depending on the fixed point $x_0$ are in a sense worse than those that depend on the variable
point $x$ since the former must be evaluated in $L^\infty(\RR_x)$ norm. To overcome this apparent
obstacle we average over $S^1$ (namely, $x_0$ in the range $(x-1,x)$) to obtain
\begin{equation}
\begin{array}{lll}
&\N_te_{\b j}(t,x)
=\cr\cr
&\q\;\dis
\sum_{q=1}^{2n}\bigg(\int_{S^1}h(\N_te_{\b j},e_q)(t,x_0)dx_0\Big) \; e_q(t,x)
\cr\cr
&\q\;\dis+\Big[\sum_{\a,\be}\Big(\big[a^\a a^{\b\be} R_{\a\be{\b j}}^{\phantom{\a\be\b j}q}+
\half(a^{\b\a}a^{\b\be}+a^{\a}a^{\be}) R_{\b\a \be{\b j}}^{\phantom{\a\be\b j}q}
\big](t,x)
\cr\cr
&\qq\;\qq\q\dis-
\int_{S^1}\big[a^\a a^{\b\be} R_{\a\be{\b j}}^{\phantom{\a\be\b j}q}+
\half(a^{\b\a}a^{\b\be}+a^{\a}a^{\be}) R_{\b\a \be{\b j}}^{\phantom{\a\be\b j}q}
\big](t,x_0)dx_0\Big)
\cr\cr
&\q\q\q\dis-\int_{[x_0,x]}\sum_{\a,\be}\big[a^\a a^{\b\be}a^s
R_{\a\be{\b j}}^{\phantom{\a\be\b j}q}{}_{,s}+
\half(a^{\b\a}a^{\b\be}+a^{\a}a^{\be}a^s) R_{\b\a \be{\b j}}^{\phantom{\a\be\b j}q}{}_{,s}
\big](t,y)dy\bigg]e_q(t,x)
\cr\cr
&
\q\qq\qquad\qquad\;\;\!
=:\big(T_{\b j}^{q}
+P_{\b j}^{q}-\int_{S^1}P_{\b j}^{q}(t,x_0)dx_0+Q_{\b j}^q\big)e_q.\cr
\end{array}
\end{equation}
Note that according to (\ref{ChristoffelTermEq}) we have
$\N_t e_{\b j}=a^{\b k}_{,x}\Gamma^p_{k\b j} e_p$, hence
\begin{equation}
\label{ChristoffelThirdEq}
h(\N_te_{\b j},e_q)=a^{\b k}_{,x}\Gamma^q_{k\b j}.
\end{equation}
Note that in (\ref{ChristoffelThirdEq}) the left hand side, hence also the right
hand side, are bona fide functions on $S^1$ (even though each term separately in the product on the right hand side
is not). Therefore,
$$
\int_{S^1}h(\N_te_{\b j},e_q)(t,x_0)dx_0
=
\int_{S^1} a^{\b k}_{,x}\Gamma^q_{k\b j}dx_0
=
-\int_{S^1} a^{\b k}\Gamma^q_{k\b j,x}dx_0
=
-\int_{S^1} a^{\b k}a^p\Gamma^q_{k\b j,p}dx_0.
$$

Switching to complex notation, as in (\ref{NLSEq}), we have
\begin{equation}
\label{NLSCircleOneEq}
\i\bfPhi_{,t}=\bfPhi_{,xx}-{\bf Q}\cdot\bfPhi-\bfS\cdot\bfPhi+{\bf W}\cdot\bfPhi-\bfT\cdot\bfPhi,
\end{equation}
\begin{equation}
\label{NLSCircleTwoEq}
\bfPhi^\a(t,x+l)=P(u(t)|_{[0,1]})^{-l}\bfPhi^\a(t,x),
\end{equation}
where ${\bf Q}:=\iota^{-1}(T_{\b j}^{q}), \; {\bf S}:=\iota^{-1}(P_{\b j}^{q}),\;
{\bf T}:=\iota^{-1}(Q_{\b j}^{q}),\; {\bf W}:=\iota^{-1}(\int_{S^1}P_{\b j}^{q}(t,x_0)dx_0)$,
and
\begin{equation}
\label{EstimatesOnNonlinearTermsEq}
\begin{array}{lll}
||{\bf Q}||<C\int_{S^1} ||\bfPhi||^2dx, \q
||{\bf W}||<C\int_{S^1} ||\bfPhi||^2dx,
\cr\cr
||{\bf S}||<C ||\bfPhi||^2dx, \q
||{\bf T}||<C\int_{S^1} ||\bfPhi||^3dx,
\end{array}
\end{equation}
where $C>0$ depends only on the geometry of $(N,h)$. We emphasize
that Equations (\ref{NLSCircleOneEq})-(\ref{NLSCircleTwoEq}) are
on $I\times \RR$. Put
$$
P(u(t)|_{[0,1]})=:e^{\i\th(t)}\in U(1), \qq\theta\in\RR,
$$
and set
$$
a:=\bfPhi^1=a^1+\i a^{\bar 1}
$$
(note that we cannot restrict $\theta$ to $[-\pi,\pi)$ in order not to violate continuity of $\theta$).
Note that following an earlier remark $\bfQ=Q^1_1, \bfS=S^1_1,
{\bf W}=W^1_1, \bfT=T_1^1$ are $\Z$-invariant. Also
\begin{equation}
\label{FirstChangeVariablesEq}
\vp(t,x):=e^{\i\th x}a(t,x)
\end{equation}
is $\Z$-invariant.    
Moreover, so are all of its $x$-derivatives. To wit,
\begin{equation}
\label{PeriodicFnEq}
\vp(t,x)_{,x}=\i\th\vp(t,x)+e^{\i\th x}a(t,x)_{,x}=\vp(t,x+1)_{,x}
\end{equation}
since $(e^{\i\th}a(t,x+1))_{,x}=a(t,x)_{,x}$, and the claim now follows by induction.
It follows that the estimates we will obtain for $\vp$ will imply the same estimates for $a$.

After the change of variable (\ref{FirstChangeVariablesEq}) Equation
(\ref{NLSCircleOneEq}) becomes
\begin{equation}
\i\vp_{,t}=\vp_{,xx}-2\i\th\vp_{,x}-(\th^2+x\th_{,t}+Q^1_1+S_1^1-W^1_1+T_1^1)\vp.
\end{equation}
Let $\beta:I\times\R\ra I\times\R$ be given by
$$
\beta(t,x)=(t,x-2\int_{[0,t]}\th ds).
$$
Let
$$
\tilde x:=x+2\int_{[0,t]}\th ds.
$$
Writing
$$
(t,x)=\beta(t,x+2\int_{[0,t]}\th ds)=\beta(t,\tilde x),
$$
Equation (\ref{NLSCircleOneEq}) becomes
\begin{equation}
\label{SoneSchrodingerEq}
\begin{array}{lll}
\i(\vp\circ\beta)_{,t}(t,\tilde x)
& = \;
         (\vp\circ\beta)_{,\tilde x\tilde x}(t,\tilde x)
-\Big(\th^2(t)+(\tilde x-{\textstyle2\int_{[0,t]}}\th ds)\th_{,t}(t)
\cr \cr
&\qq+\,(Q_1^1\circ\beta+S_1^1\circ\beta-W_1^1\circ\beta+T_1^1\circ\beta)(t,\tilde x)\Big)
(\vp\circ\beta)(t,\tilde x).
\end{array}
\end{equation}
This equation is on $I\times S^1$.

\begin{remark}
{\rm
Note that here it was crucial that $\th$ does not depend on $x$ in order to have $\frac{\partial \tilde x}{\partial x}=1$.
This is also the difference from the situation in Equation (\ref{NLSVariableMetricEq}).
}
\end{remark}

The main result of this section is:

\begin{theo}
\label{SThm}
Let $(M,g)=(S^1,dx\otimes dx)$ and let $(N,J,h)$ be a complete Riemann surface with bounded geometry.
Then the system of NLS equations (\ref{SoneSchrodingerEq})
is locally well-posed in the space $L^4(\RR,L^4(S^1,\RR^{2n}))$.
\end{theo}

This will then imply:
\begin{cor}
\label{SCor}
Let $(M,g)=(S^1,dx\otimes dx)$ and let $(N,J,h)$ be
a complete Riemann surface with bounded geometry. 
Then for integers $k\ge 2$ 
equation (\ref{SchEq}) with $u_0\in W^{k,2}(S^1,N)$ admits a
unique solution $u\in C^0(\RR,W^{k,2}(S^1,N))$. In particular $u$
is smooth if $u_0$ is.
\end{cor}

We turn to the proof of Theorem \ref{SThm}.

\begin{proof}
We will use Equation (\ref{SoneSchrodingerEq}) in order to obtain a priori estimates on
$$
\tilde\vp(t,\tilde x):=\vp\circ\beta(t,\tilde x).
$$
The estimates on $\tilde\vp$ and on $\vp$ are equivalent since the two functions only differ by
a time-dependent translation in the space direction.
We will localize in time: indeed it is enough to prove existence of local (in time) solutions
of Equation (\ref{SoneSchrodingerEq}) in
$C^0(\RR_{t,loc},L^2(S^1))\cap L^4(\RR_{t,loc}\times S^1)$
depending in a good manner only on $||\tilde\vp||_{L^2(\RR_{\tilde x})}=||a||_{L^2(\RR_x)}$ and
a bounded constant depending on time, since that
will rule out finite-time blow-up.

Let us now recall some of the work of Bourgain that will be of central importance in what follows
\cite{B} (see also \cite{G}).
First, we recall the following Fourier restriction estimates of Bourgain:

\begin{lem}
\label{FirstFourierAnalysisLemma}
{\rm (See \cite{B}, p. 112.)}
Let $\vp$ be a periodic solution of the linear \Sch equation on $S^1$.
Then
$$
||\vp||_{L^4(S^1\times S^1)}\le \sqrt 2||\vp(0)||_{L^2(S^1)},
$$
and dually
$$
||\vp||_{L^2(S^1\times S^1)}\le \sqrt 2||\vp||_{L^{4/3}(S^1\times S^1)}.
$$
\end{lem}

More generally, Bourgain proved the following fundamental result that allows for the same estimate---now
with appropriate weights---even for an arbitrary function whose Fourier modes are not necessarily
restricted to the parabola $\{(n,n^2)\,:\, n\in\ZZ\}$.
We state the result although we will only directly use a consequence of it.

\begin{lem}
\label{SecondFourierAnalysisLemma}
{\rm (See \cite{B}, Proposition 2.33.)}
Let $f(x,t)=\sum_{m,n\in\ZZ}a_{m,n}e^{\i(mx+nt)}$ be a function on $S^1\times S^1$. Then
$$
\Big(\sum_{m,n\in\ZZ}(|n-m^2|+1)^{-3/4} |a_{m,n}|^2\Big)^{1/2}\le c||f||_{L^{4/3}(S^1\times S^1)}.
$$
In addition, if $|\lambda_{m,n}|\le (1+|n-m^2|)^{-3/4}$, then
$$
||\sum_{m,n\in\ZZ}\lambda_{m,n}a_{m,n}e^{\i(mx+nt)}||_{L^4(S^1\times S^1)}\le c ||f||_{L^{4/3}(S^1\times S^1)}.
$$
In both estimates $c>0$ is some universal constant.
\end{lem}
Using this estimate Bourgain obtains the following $L^4$ estimate for the nonlinear contribution
in Duhamel's formula. This estimate will play a central r\^ole below.

\begin{lem}
\label{ThirdFourierAnalysisLemma}
{\rm (See \cite{B}, \S4)}
Let $F\in L^{4/3}(S^1\times S^1)$.
For any $0<\delta<1/8$ and $\dis0<B<\frac1{100\delta}$ there holds
$$
\Big|\Big|
\int_0^{2\delta} S(t-\tau)F(\tau,x)d\tau
\Big|\Big|_{L^4(S^1\times S^1)}
\le
C (B^{-1/4}+\delta B)||F||_{L^{4/3}(S^1\times S^1)},
$$
where $C>0$ is some universal constant.
\end{lem}
The constant $B$
can be thought of as a Fourier mode cut-off parameter, measuring distance of a lattice point in $\ZZ^2$
from the parabola $\{(m,m^2)\,:\,m\in\ZZ\}$. The constant $\delta$ is the time cut-off parameter.

Equation (\ref{SoneSchrodingerEq}) is equivalent to the integral equation
\begin{equation}
\label{IntegralSoneNLSEq}
\tilde\vp(t,\tilde x)=
S(t)\tilde\vp(0,\tilde x)
-
\i\int_0^t S(t-\tau) F(\tau,\tilde x)d\tau.
\end{equation}
with
\begin{equation}
\label{NonlinearTermEq}
\begin{array}{lll}
F(\tau,\tilde x) & =
-\Big(\th^2(\tau)+(\tilde x-{\textstyle2\int_{[0,\tau]}}\th ds)\th_{,t}(\tau)
\cr\cr
&\qq\q\;+\,(Q_1^1\circ\beta+S_1^1\circ\beta-W_1^1\circ\beta+T_1^1\circ\beta)(\tau,\tilde x)\Big)
\tilde\vp(\tau,\tilde x).
\end{array}
\end{equation}
There is a subtlety here: the time derivative of $\tilde\vp$ (or
of $\vp$) is not necessarily $\ZZ$-invariant (in $\tilde x$).
However, Equation (\ref{SoneSchrodingerEq}) holds on $I\times S^1$
and it is equivalent to the integral equation
(\ref{IntegralSoneNLSEq}).

We would like to obtain an a priori $L^4$ estimate on $\tilde\vp$. We localize in time, namely multiply
Equation (\ref{IntegralSoneNLSEq}) by a smooth cut-off function in time $\psi(t)$ satisfying
$\psi=1$ on $[-\delta,\delta]$ and $\psi=0$ for $|t|\ge 2\delta$.
Here $\delta$ is a positive number smaller than $1/8$ to be specified
later. We may thus regard $\psi\tilde\vp$ as a function on $S^1\times S^1$ with period $1$ in both the $t$ and
$\tilde x$ variables and Bourgain's estimates apply.

First, the linear term satisfies
$$
||\psi S(t)\tilde\vp(0,\tilde x)||_{L^4(S^1\times S^1)}\le \sqrt 2 ||\tilde\vp(0,\cdot)||_{L^2(S^1)}
=\sqrt{2 E(u_0)},
$$
according to Lemma \ref{FirstFourierAnalysisLemma}.

Next, we estimate the integral term. The terms involving $\bf Q$ and $\bf W$ are simpler since
$||\bf Q||$ and $||\bf W||$ are uniformly bounded according to (\ref{EstimatesOnNonlinearTermsEq}) and conservation of energy.

We now turn to the other terms. First, using Lemma \ref{FirstFourierAnalysisLemma} under
the integral sign, and assuming $||\theta||_{L^{\infty}}\le C$, we have
\begin{equation}
\label{FirstNonLinearEstimateEq}
\begin{array}{lll}
\dis||\psi\int_0^t S(t-\tau) \big(\th^2(\tau)\tilde\vp(\tau,\tilde x)\big)d\tau||_{L^4(S^1\times S^1)}
& \dis\le \int_0^{2\delta} ||\th^2(\tau)\tilde\vp(\tau,\,\cdot\,)||_{L^2(S^1)}d\tau
\cr\cr
&\dis\le 2C^2\delta||\vp(0)||_{L^2(S^1)}.
\end{array}
\end{equation}
To show that this assumption holds use the representation of the holonomy given by
(\ref{FirstHolonomyReprEq}): $|\theta(t)|\le\int_0^1|\T kij||a^i|dx\le C' E(u_0)^{1/2}$
where we have used the assumption of bounded geometry---indeed it implies that Christoffel
symbols are uniformly bounded \cite{E}.

Second, $|\tilde x|\le1$ and so $|\tilde x-{\textstyle2\int_{[0,\tau]}}\th ds|\le 1+2\cdot1\cdot C$.
Let $\{\alpha_1,\alpha_{\b1}\}$ be an orthonormal coframe dual to $\{e_1,e_{\b1}\}$.
To compute the time derivative of $\th$, recall that by Equation (\ref{HolonomyEq}) we have
$$
\th(t)=\int_{D_u}KdV_{N,h}=\int_{D_u}K \alpha_1\w\alpha_{\b 1}
=\int_{I\times S^1}K\circ u(t,x)[a^1b^{\b 1}-a^{\b 1}b^1]dx\w dt,
$$
since $u^\star \alpha_1=a^1dx+b^1dt,\; u^\star \alpha_{\b 1}=a^{\b 1}dx+b^{\b 1}dt$.
Combining this with the fact that by (\ref{FlowEqInCoordinatesEq}) we have $b^k=a^{\b k}_{,x}$,
we then have
$$
\th_{,t}=\int_{S^1} K\circ u(t,x)(a^1b^{\b 1}-a^{\b 1}b^1)dx
=-\frac12\int_{S^1}K\circ u(t,x)((a^1)^2+(a^{\b1})^2)_{,x}dx.
$$
Integrating by parts this becomes
$$
\th_{,t}=\frac12\int_{S^1}(K\circ u(t,x))_{,x}((a^1)^2+(a^{\b1})^2)dx
=\frac12\int_{S^1}K_{,s}\circ u(t,x)a^s((a^1)^2+(a^{\b1})^2)dx.
$$
By bounded geometry we therefore have
\begin{equation}
\label{SecondSecondNonLinearEstimateEq} ||\th_{,t}||_{L^\infty}\le
C||a||^3_{L^3(S^1)}.
\end{equation}
Therefore the term $(\tilde x-{\textstyle2
\int_{[0,\tau]}}\th ds)\th_{,t}(\tau) \tilde\vp(\tau,\tilde x)$ behaves in
the same way as the term
$T_1^1\circ\beta\tilde\vp(\tau,\tilde x)$ in (\ref{NonlinearTermEq}), and so
it's enough to treat the latter. We will do that shortly.

Third, $|S_1^1\circ\beta\cdot \tilde\vp|<C
|\tilde\vp|^3$, and therefore this term may be estimated in
$L^4(S^1\times S^1)$ just like in Bourgain's estimates for a cubic
nonlinearity. More precisely, by Lemma
\ref{ThirdFourierAnalysisLemma} we have
\begin{equation}
\label{ThirdNonLinearEstimateEq}
\dis\Big|\Big|\psi\int_0^t S(t-\tau)
\big(|\tilde\vp|^3(\tau,\tilde x)\big)d\tau
\Big|\Big|_{L^4(S^1\times S^1)}
\le
C
(\delta B+B^{-1/4}) ||\psi\tilde\vp||^3_{L^4(S^1\times S^1)}
,
\end{equation}
with $B>0$ as in the Lemma, $\delta$ is as before the time cut-off parameter,
and $C>0$ is a uniform constant.

Fourth, using Lemma \ref{FirstFourierAnalysisLemma} and energy conservation we have
\begin{equation}
\label{FourthNonLinearEstimateEq}
\begin{array}{lll}
&\dis\Big|\Big|\psi\int_0^t S(t-\tau)
\big(\tilde\vp\int_{S^1}|\tilde\vp|^3d\tilde x(\tau,\tilde x)\big)d\tau\Big|\Big|_{L^4(S^1\times S^1)}
\cr\cr
& \qquad\qquad
\dis\le C\int_0^{2\delta}\Big|\Big|
\tilde\vp(\tau,\,\cdot\,) \int_{S^1}|\tilde\vp(\tau,\,\cdot\,)|^3d\tilde x
\Big|\Big|_{L^2(S^1)}d\tau
\cr\cr
&\dis \qquad\qquad
\le C' \int_0^{2\delta}\int_{S^1} |\tilde\vp(\tau,\,\cdot\,)|^3d\tilde x d\tau
\le C' \root 4 \of 2 \delta^{1/4} ||\vp||_{L^4(S^1\times S^1)}^3.
\end{array}
\end{equation}

Combining Equations (\ref{FirstNonLinearEstimateEq})--%
(\ref{FourthNonLinearEstimateEq}) we have
\begin{equation}
||\psi\tilde\vp||_{L^4(S^1\times S^1)}
\le
C\big(
(1+\delta)||\vp||_{L^2(\RR_x)}+(\delta B+B^{-1/4}+\delta^{1/4}) ||\psi\tilde\vp||^3_{L^4(S^1\times S^1)}
\big).
\end{equation}
By choosing $B$ large enough and then choosing $\delta$ small enough, in such a manner that $\delta B$ is
also small enough (all these choices depend only on the initial energy and the background geometry)
we therefore may argue as Bourgain (a standard Picard iteration argument, see \cite{B}, p. 139) to obtain
a uniform estimate on $||\tilde\vp||_{L^4(\RR_{t,loc}\times S^1)}=||\vp||_{L^4(\RR_{t,loc}\times S^1)}
=||a||_{L^4([-\delta,\delta]\times S^1)}$.
Global well-posedness (existence and uniqueness) in $L^4(\RR_{t,loc}\times S^1)$ now follows from energy conservation.

To obtain global well-posedness for our original flow equation in $W^{k,2}$
we take $k-1$ derivatives of Equation
(\ref{SoneSchrodingerEq}).
Unlike in the case of $M=\RR$ (Section \ref{SectionR}), for some of the terms
we might need to use some interpolation inequalities already for the case $k=2$.
In addition, compared to Bourgain's situation of the cubic NLS on $S^1$ we obtain in fact
certain terms that are worse than those obtained by differentiating the cubic NLS.
For example, for $k=2$ we obtain terms of order
$|\tilde\vp |^4$ and $|\tilde\vp_{,x}\tilde\vp^3|$.
Such terms may be handled nevertheless. We carry this in some detail in the case $k=2$, omitting
the details in the case $k\ge3$ as they are similar.

Set $c:=\tilde\vp_{,\tilde x}$. Taking a derivative of Equation (\ref{SoneSchrodingerEq})
we obtain
$$
\i c_{,t}=c_{,\tilde x\tilde x}+F
$$
where
\begin{equation}
\label{TermsInNonlinearTermOfDerivativeEq}
|F|<C\big(|\tilde\vp|^4+|c|(1+|\tilde\vp|^2+|\tilde\vp|^3+||\tilde\vp||_{L^1(S^1)}^3) \big).
\end{equation}
As before we would like to obtain an $L^4(\RR_{t,loc}\times S^1)$ estimate, this time
for $c$. If we attempt to use Lemma \ref{FirstFourierAnalysisLemma} in order to handle the term $|\vp|^4$ we
would need to estimate
$$
\int_0^t ||\vp(\tau,\,\cdot\,)^4 ||_{L^2(S^1)} d\tau.
$$
Using the Gagliardo-Nirenberg inequality
\cite{A}, p. 93, we would have
$$
||\vp||^4_{L^8(S^1)}\le c ||\nabla_x\vp||_{L^4(S^1)}^{1/2}||\vp||_{L^4(S^1)}^{7/2}
$$
and this can be estimated by
$
\epsilon (||\nabla_x\vp||_{L^4(S^1)}+1)^{4}+C_\epsilon (||\vp||_{L^4(S^1)}+1)^{4},
$
for appropriate $\epsilon,C_\epsilon>0$.
This is not enough for our purposes since we need to control the coefficient multiplying
the second term. Instead, we apply Lemma \ref{ThirdFourierAnalysisLemma} which gives an improved
estimate. Namely, we obtain
\begin{equation}
\label{FirstFirstTermInFEq}
\Big|\Big|\psi\int_0^t S(t-\tau)
\big(|\tilde\vp|^4(\tau,\tilde x)\big)d\tau
\Big|\Big|_{L^4(S^1\times S^1)}
\le
C(B^{-1/4}+\delta B)||\vp||^4_{L^{16/3}(\RR_{t,loc}\times S^1)}.
\end{equation}
By the Gagliardo-Nirenberg inequality, as before, this is therefore bounded by
\begin{equation}
\label{FirstSecondTermInFEq}
C(B^{-1/4}+\delta B)\Big(
\epsilon ||\nabla_x\vp||_{L^4(R_{t,loc}\times S^1)}^{4}+C_{\epsilon} ||\vp||_{L^4(R_{t,loc}\times S^1)}^{4}.
\Big)
\end{equation}
Next, to handle the term $|c||\tilde\vp|^3\le \epsilon|c|^2+C_\epsilon|\vp|^6$,
again by Lemma \ref{ThirdFourierAnalysisLemma}, we compute
\begin{equation}
\label{SecondTermInFEq}
\begin{array}{lll}
&\Big|\Big|\psi\int_0^t S(t-\tau)
\big(|c\tilde\vp|^3(\tau,\tilde x)\big)d\tau
\Big|\Big|_{L^4(S^1\times S^1)}
\cr\cr
&\le
C(B^{-1/4}+\delta B)
\Big(
C_{\epsilon'}||\vp||^6_{L^8(\RR_{t,loc}\times S^1)}
+
\epsilon'||c||^2_{L^{8/3}(\RR_{t,loc}\times S^1)}
\Big)
\cr\cr
&\le
C(B^{-1/4}+\delta B)
\Big(
C_{\epsilon'}C_{\epsilon''}||\vp||^6_{L^4(\RR_{t,loc}\times S^1)}
+
C_{\epsilon'}\epsilon''||c||^6_{L^4(\RR_{t,loc}\times S^1)}
\cr\cr
&\qquad\qquad\qquad\qquad\qquad
+
\epsilon'C_{\epsilon'''}||\vp||^2_{L^4(\RR_{t,loc}\times S^1)}
+
\epsilon'\epsilon'''||c||^2_{L^4(\RR_{t,loc}\times S^1)}
\Big).
\end{array}
\end{equation}

The other terms in (\ref{TermsInNonlinearTermOfDerivativeEq}) are of lower order than the two terms we just estimated.
Since we are free to choose appropriate $B$ and $\delta$,
it follows from (\ref{FirstFirstTermInFEq})--(\ref{SecondTermInFEq}) and the Duhamel formula that
the same Picard iteration type argument that was used earlier applies to our situation.
Note that this estimate depends on the $L^2(S^1)$ norm of $c$, which
we control under our assumption $\tilde\vp(0)\in W^{1,2}(S^1)$ and on the $L^4(\RR_{t,loc}\times S^1)$
norm of $\vp$ that we already control uniformly.
This concludes the proof of Theorem \ref{SThm}.
\end{proof}

\section{Maps of the circle into a \K manifold}
\label{SectionCircleHigherD}

In this section we explain the difficulties that one encounters if one tries to apply the same methods
to treat the case of maps from the circle to K\"ahler manifolds of arbitrary dimension  $n\ge 1$.

For general $n$, one gets an expression for a solution of
(\ref{GaugeEq}) given by the chronological exponential
\begin{equation}
\label{AEq}
\begin{array}{lll}
e(t,x+1) & =A^{-1}(t,x)e(t,x)
\cr\cr
& :=\lim_{n\ra\infty}
e^{-\frac1{n}B_U(t,x+1)}e^{-\frac1{n}B_U(t,x+\frac{n-1}n)}\cdots
e^{-\frac1{n}B_U(t,x+\frac1n)}e(t,0)
\end{array}
\end{equation}
(see for example \cite{DFN}).

First we observe that $A(t,x)$ does not depend on $x$: indeed applying $\N_x$ to Equation (\ref{AEq})
and using the fact that $\nabla_x e=0$ we obtain $A_{,x}=0$. From now on we simply write $A(t)$.

Now $\bfPhi(t,x+1)=A(t)\bfPhi(t,x)$. Since $A(t)$ is unitary (and hence normal) it is unitarily diagonalizable and we set
$$
A(t)=U(t)^\star D(t)U(t),
$$
with $U(t)\in U(n)$ and $D(t)=\diag(e^{\i\th_1},\ldots,e^{\i\th_n})$.
Let
$$
\bfPhitilde(t,x):=U(t)^\star D(t)^{-x} U(t)\bfPhi(t,x)=A(t)^{-x}\bfPhi(t,x).
$$
The vector-valued function $\bfPhitilde$ is periodic in $x$. Moreover, by a computation similar to
(\ref{PeriodicFnEq}), so are all of its $x$-derivatives. We have
\begin{equation}
\begin{array}{lll}
\label{bfPhitDerivEq}
\bfPhi_{,t} & =(A(t)^x)_{,t}\bfPhitilde+A(t)^x\bfPhitilde_{,t},
\cr\cr
\bfPhi_{,x} & =U(t)^\star D(t)^x \diag(\i\th_i)U(t)\bfPhitilde+U(t)^\star D(t)^xU(t)\bfPhitilde_{,x},
\cr\cr
\bfPhi_{,xx} & =U(t)^\star D(t)^x \diag(-\th^2_i)U(t)\bfPhitilde+2U(t)^\star D(t)^x\diag(\i\th_i)U(t)\bfPhitilde_{,x}
\cr\cr
&\q +U(t)^\star D(t)^xU(t)\bfPhitilde_{,xx}.
\end{array}
\end{equation}
It follows that
\begin{equation}
\begin{array}{lll}
\i\bfPhi_{,t}-\bfPhi_{,xx}
& =A(t)^x
\Big[
\i\bfPhitilde_{,t}
-
\bfPhitilde_{,xx}+(A(t)^x)_{,t}\bfPhitilde
\cr\cr
&\q-U(t)^\star\diag(-\th^2_i)U(t)\bfPhitilde-2U(t)^\star\diag(\i\th_i)U(t)\bfPhitilde_{,x}
\Big]
\end{array}
\end{equation}
Therefore equations (\ref{NLSCircleOneEq})-(\ref{NLSCircleTwoEq}) may be rewritten as
\begin{equation}
\begin{array}{lll}
\i\bfPhitilde_{,t}
& =
\bfPhitilde_{,xx}
-(A(t)^x)_{,t}\bfPhitilde+U(t)^\star\diag(-\th^2_i)U(t)\bfPhitilde
\cr\cr
&\q+2U(t)^\star\diag(\i\th_i)U(t)\bfPhitilde_{,x}
-A(t)^{-x}\big({\bf Q}\cdot\bfPhi+\bfS\cdot\bfPhi-{\bf W}\cdot\bfPhi+\bfT\cdot\bfPhi\big).
\end{array}
\end{equation}
Note that the last term is expressed in terms of $\bfPhi$ instead of $\tilde\bfPhi$. However as far
as the estimates are concerned this is not important since it
involves no derivatives and the two vectors differ by a unitary transformation.
Two problems now arise. First, one needs to obtain an estimate on the variation of the holonomy matrix $A(t)$
along the flow. Such an estimate was available in the one-dimensional setting due to the Gauss-Bonnet theorem.
Second,
the matrix multiplying the first derivative term is not diagonal and so it is not clear how to eliminate this term.

Although this requires some work, and we will not attempt to
provide the details here, the first difficulty may be overcome
using the theory developed by Chacon and Fomenko for a
non-commutative version of the Stokes' Theorem product integrals
\cite{CF} (see also the classical references \cite{N,S}). To
approach the second difficulty one may consider $\hat
\bfPhi:=U(t)^\star D(t)^{-x}\bfPhi$ instead of $\tilde\bfPhi$.
Then the matrix multiplying the first derivative of $\hat \bfPhi$
is diagonal. Therefore, we may apply the space-time transformation
as in the Riemann surface case, however for each equation in the
system separately. However, this introduces a new obstacle.
Indeed, then one needs to control the time derivative of $D(t)$ as
well as of $U(t)$. The main difficulty comes from the latter. In
general, the unitary diagonalizing matrix does not vary smoothly
(or even continuously) even when a family of matrices does (see
\cite{K}, p. 111). Instead one may try to diagonalize $A(t)$
smoothly.
However, to the best of our knowledge, even given such a
diagonalization, the problem is that even though the diagonalizing
matrix is then smooth one has essentially no control over its
derivatives (i.e., estimates on these derivatives in terms of
derivatives of $A(t)$). We hope to come back to this problem in
the future. In some sense the two transformations (to $\tilde
\bfPhi$ and to $\hat\bfPhi$) are dual to each other, and one may
ask whether for higher-dimensional domains the two troublesome
terms, namely the first derivative term and the derivative of the
holonomy, may be a source for finite-time blow-up.

\end{document}